\documentclass[11pt,a4paper,twoside]{article}%
\usepackage{amsfonts}
\usepackage{amsmath}
\usepackage{geometry}
\usepackage{version}
\usepackage{amssymb}
\usepackage{graphicx}%
\providecommand{\U}[1]{\protect\rule{.1in}{.1in}}
\newtheorem{theorem}{Theorem}
\newtheorem{acknowledgement}[theorem]{Acknowledgement}

\newtheorem{definition}[theorem]{Definition}
\newtheorem{example}[theorem]{Example}

\newenvironment{proof}[1][Proof]{\noindent\textbf{#1.} }{\ \rule{0.5em}{0.5em}}
\begin{document}

\title{Bispectrum for non-Gaussian homogeneous and isotropic field on the plane}
\author{Gy\"{o}rgy Terdik\\Faculty of Informatics, University of Debrecen, Hungary, \\Terdik.Gyorgy$@$inf.unideb.hu}
\maketitle

\begin{abstract}
The object of this paper is to characterize the third order moments
(cumulants) and bispectra of a homogeneous isotropic field defined on a plane.
We establish a one to one correspondence between the third order cumulants and
the bispectra of such a process in terms of Bessel functions.

\end{abstract}

Bispectrum

\begin{center}
\bigskip\textit{Dedicated to Lajos Tam\'{a}ssy on the occasion of his 90th
birthday.}
\end{center}

\section{Introduction}

In many real applications associated with random fields, the assumption of
Gaussianity may be sometimes unrealistic. For example,consider the data of
cosmic microwave background (CMB) anisotropies provided by NASA which some
scientists believe to be non-Gaussian, \cite{Marinucci2011}. Although CMB data
given are on a surface of the sphere there are problems concerning on the
primordial field on the whole space including the investigation of the
bispectrum as well, see \cite{verde2000large}, \cite{yadav2007fast},
\cite{achitouv2012primordial}. In time series analysis the non-Gaussianity has
been well studied \cite{SRao-Gabr-test-80}, \cite{Rao84a},
\cite{Hinich_test-82}, \cite{Terdik-Math-JTSA-98}, \cite{TerdikLN99}. It is
known that for a Gaussian time series the bispectrum and all higher order
spectra greater then two are zero, and equally well known is the fact that for
a non Gaussian process defined on a real line the higher order cumulant
spectra and higher order cumulants uniquely determine each other
\cite{Brill-book-01}. No such results are available and well known for a
homogeneous isotropic field defined on a plane.

In this paper we consider homogeneous and isotropic fields which are not
necessarily Gaussian. The second order properties of such a field are well
known, both the covariance and the spectrum depend on the distance between the
locations and the wave numbers respectively, \cite{Whittle1954},
\cite{Yadrenko1983}. The third order covariances (third order cumulants)
depends on three locations and because of the invariance under   shifting
(homogeneity) and invariance under the rotation (isotropy), it depends only on
the distances between locations. In other words the rigid body movement keeps
the triangle defined by the three locations fixed. We show that the third
order covariances define the bispectrum which depends on three wave numbers
forming a triangle.  The main result of this paper describes the unique
relation between the third order covariances and the bispectrum of a
homogeneous and isotropic field on the plane.

\subsection{Homogeneous and isotropic field on the plane}

We consider a homogeneous real valued stochastic field $X\left(
\underline{x}\right)  $ on $\mathbb{R}^{2}$ with $EX\left(  \underline{x}%
\right)  =0$. Let us suppose that $X\left(  \underline{x}\right)  $ is
continuous (in mean square sense), its spectral representation is
\begin{equation}
X\left(  \underline{x}\right)  =\int_{\mathbb{R}^{2}}e^{i\underline{x}%
\cdot\underline{\omega}}Z\left(  d\underline{\omega}\right)  ,\quad
\underline{\omega},\underline{x}\in\mathbb{R}^{2},\label{X_hom}%
\end{equation}
with a finite spectral measure%
\[
E\left\vert Z\left(  d\underline{\omega}\right)  \right\vert ^{2}=F_{0}\left(
d\underline{\omega}\right)  .
\]
By homogeneity we mean (in strict sense ) the distribution of $X\left(
\underline{x}\right)  $ is translation invariant, see \cite{Yaglom-book-87}
for details. Rewrite $X\left(  \underline{x}\right)  $ in terms of polar
coordinates
\[
X\left(  r,\varphi\right)  =\int_{0}^{\infty}\int_{0}^{2\pi}e^{i\rho
r\cos\left(  \varphi-\eta\right)  }Z\left(  \rho d\rho d\eta\right)  ,
\]
where $\underline{x}=\left(  r,\varphi\right)  $, $\underline{\omega}=\left(
\rho,\eta\right)  $ are polar coordinates, $r=\left\vert \underline{x}%
\right\vert =\sqrt{x_{1}^{2}+x_{2}^{2}}$, and $\rho=\left\vert
\underline{\omega}\right\vert $, $\underline{x}\cdot\underline{\omega}%
=r\rho\cos\left(  \varphi-\eta\right)  $. Now we use the Jacobi-Anger
expansion, see \cite{Arfken2001} Sect 11 ,
\begin{equation}
e^{i\rho r\cos\left(  \varphi-\eta\right)  }=\sum_{\ell=-\infty}^{\infty
}i^{\ell}J_{\ell}\left(  \rho r\right)  e^{i\ell\left(  \varphi-\eta\right)
},\label{Expans_Jacobi_Anger}%
\end{equation}
and substitute it into the above spectral representation of $X\left(
r,\varphi\right)  $
\begin{equation}
X\left(  r,\varphi\right)  =\sum_{\ell=-\infty}^{\infty}e^{i\ell\varphi}%
\int_{0}^{\infty}J_{\ell}\left(  \rho r\right)  Z_{\ell}\left(  \rho
d\rho\right)  \label{Repr_spherical}%
\end{equation}
where $J_{\ell}$ denotes the Bessel function of the first kind
\cite{Abramowit12}, and the series of stochastic spectral measures $Z_{\ell
}\left(  \rho d\rho\right)  $ is connected to the $Z\left(  d\underline{\omega
}\right)  $ by the integral
\[
Z_{\ell}\left(  \rho d\rho\right)  =\int_{0}^{2\pi}i^{\ell}e^{-i\ell\eta
\ }Z\left(  \rho d\rho d\eta\right)  .
\]
The representation (\ref{Repr_spherical}) will be an orthogonal (uncorrelated)
representation if we assume that $F_{0}\left(  d\underline{\omega}\right)  $
is isotropic (invariant under rotations) $F_{0}\left(  d\underline{\omega
}\right)  =E\left\vert Z\left(  d\underline{\omega}\right)  \right\vert
^{2}=E\left\vert Z\left(  \rho d\rho d\eta\right)  \right\vert ^{2}=F\left(
\rho d\rho\right)  d\eta$. Note here that the general theorem of Yadrenko
(\cite{Yadrenko1983} Theorem 1. pp.5) on the spectral representation of a
homogeneous and isotropic field $X\left(  r,\varphi\right)  $ gives the
representation (\ref{Repr_spherical}) with real valued stochastic spectral
measures $Z_{\ell}\left(  \cdot\right)  $ constructed directly form the field
$X\left(  r,\varphi\right)  $ itself. The stochastic spectral measures
$Z_{\ell}\left(  \rho d\rho\right)  $ defined above is complex valued and has
the following property
\begin{equation}
Z_{-\ell}\left(  \rho d\rho\right)  =\left(  -1\right)  ^{\ell}\overline
{Z_{\ell}\left(  \rho d\rho\right)  }.\label{Z_compl}%
\end{equation}
Indeed, since the field $X\left(  \underline{x}\right)  $ is real valued we
have
\begin{align*}
X\left(  r,\varphi\right)   &  =\sum_{\ell=-\infty}^{\infty}e^{i\ell\varphi
}\int_{0}^{\infty}J_{\ell}\left(  \rho r\right)  Z_{\ell}\left(  \rho
d\rho\right)  \\
&  =\sum_{\ell=-\infty}^{\infty}e^{-i\ell\varphi}\int_{0}^{\infty}J_{\ell
}\left(  \rho r\right)  \overline{Z_{\ell}\left(  \rho d\rho\right)  }.
\end{align*}
Moreover, from the well known formula $J_{\ell}\left(  \cdot\right)  =\left(
-1\right)  ^{\ell}J_{-\ell}\left(  \cdot\right)  $ we have
\begin{align*}
\sum_{\ell=-\infty}^{\infty}e^{-i\ell\varphi}\int_{0}^{\infty}J_{\ell}\left(
\rho r\right)  \overline{Z_{\ell}\left(  \rho d\rho\right)  } &  =\sum
_{\ell=-\infty}^{\infty}e^{-i\ell\varphi}\int_{0}^{\infty}J_{-\ell}\left(
\rho r\right)  \left(  -1\right)  ^{\ell}\overline{Z_{\ell}\left(  \rho
d\rho\right)  }\\
&  =\sum_{\ell=-\infty}^{\infty}e^{-i\ell\varphi}\int_{0}^{\infty}J_{-\ell
}\left(  \rho r\right)  Z_{-\ell}\left(  \rho d\rho\right)  ,
\end{align*}
hence (\ref{Z_compl}) follows. The identity (\ref{Z_compl}) implies  that the
$\ell^{th}$ term and the $-\ell^{th}$ terms of the expansion
(\ref{Repr_spherical}) are conjugates of each other. Moreover $Z_{\ell}\left(
\cdot\right)  $ is orthogonal%
\begin{align}
\operatorname*{Cov}\left(  Z_{\ell_{1}}\left(  \rho_{1}d\rho_{1}\right)
,Z_{\ell_{2}}\left(  \rho_{2}d\rho_{2}\right)  \right)   &  =\int_{0}^{2\pi
}i^{\left(  \ell_{1}-\ell_{2}\right)  }e^{-i\left(  \ell_{2}-\ell_{1}\right)
\eta}F_{0}\left(  d\eta\rho d\rho\right)  \nonumber\\
&  =\delta_{\ell_{1}-\ell_{2}}2\pi F\left(  \rho d\rho\right)  ,\label{Z_cov}%
\end{align}
where $\delta_{\ell_{1}-\ell_{2}}$ is the Kronecker delta. Note that the
spectral measure $F\left(  \rho d\rho\right)  $ of $Z_{\ell}\left(  \rho
d\rho\right)  $ does not depend on $\ell$. We shall assume in particular cases
that $F\left(  \rho d\rho\right)  $ is absolutely continuos, i.e. $F\left(
\rho d\rho\right)  =\sigma^{2}\left\vert A\left(  \rho\right)  \right\vert
^{2}\rho d\rho$, here $\sigma^{2}\left\vert A\left(  \rho\right)  \right\vert
^{2}$ is usually known as the second order spectrum. In view of this
observation, we can rewrite $X\left(  r,\varphi\right)  $ in terms of white
noise measures $W_{\ell}\left(  \rho d\rho\right)  $ with constant spectrum%
\[
\operatorname*{Cov}\left(  W_{\ell_{1}}\left(  \rho_{1}d\rho_{1}\right)
,W_{\ell_{2}}\left(  \rho_{2}d\rho_{2}\right)  \right)  =\delta_{\ell_{1}%
-\ell_{2}}\sigma^{2}\rho d\rho,
\]
hence (\ref{Repr_spherical}) becomes
\[
X\left(  r,\varphi\right)  =\sum_{\ell=-\infty}^{\infty}e^{i\ell\varphi}%
\int_{0}^{\infty}J_{\ell}\left(  \rho r\right)  A\left(  \rho\right)  W_{\ell
}\left(  \rho d\rho\right)  .
\]

\section{Isotropy on the plane}

We consider rotations about the origin of the coordinate system. Under a
rotation (passive) $g\in SO\left(  2\right)  $, we mean a rotation when
vectors remain fixed, but the point it defines is given by a new set of
coordinates. A rotation $g$ is characterized by an angle $\gamma$ and by the
rotation matrix
\[
g=%
\begin{bmatrix}
\cos\gamma & -\sin\gamma\\
\sin\gamma & \cos\gamma
\end{bmatrix}
.
\]
If $\underline{x}\in\mathbb{R}^{2}$ is given in polar coordinates
$\underline{x}=\left(  r,\varphi\right)  $, then $g\underline{x}=\left(
r,\varphi-\gamma\right)  $, and as usual the operator $\Lambda\left(
g\right)  $ acts on functions $f\left(  r,\varphi\right)  $, such that
$\Lambda\left(  g\right)  f\left(  r,\varphi\right)  =f\left(  g^{-1}\left(
r,\varphi\right)  \right)  =f\left(  r,\varphi+\gamma\right)  $.

The isotropy  usually is defined through the invariance of the covariance
structure. This is satisfactory for Gaussian cases but for non-Gaussian fields
we need invariance of higher order cumulants as well. We use a stronger
definition to achieve a similar invariance to be able to define third order
spectrum, which we will propose below.

\begin{definition}
A homogeneous stochastic field $X\left(  \underline{x}\right)  $ is strictly
isotropic if all finite dimensional distributions of $X\left(  \underline{x}%
\right)  $ are invariant under rotation.
\end{definition}

If the homogeneous field $X\left(  \underline{x}\right)  $ is Gaussian,then
the isotropy of the spectral measure $F_{0}\left(  d\underline{\omega}\right)
$, i.e. in polar coordinates $F_{0}\left(  d\underline{\omega}\right)
=F\left(  \rho d\rho\right)  d\eta/\left(  2\pi\right)  $, implies
\[
\operatorname*{Cov}\left(  \Lambda\left(  g\right)  X\left(  \underline{x}%
_{1}\right)  ,\Lambda\left(  g\right)  X\left(  \underline{x}_{2}\right)
\right)  =\operatorname*{Cov}\left(  X\left(  \underline{x}_{1}\right)
,X\left(  \underline{x}_{2}\right)  \right)  ,
\]
for each $\underline{x}_{1}$, $\underline{x}_{2}$ and for every $g\in
SO\left(  2\right)  $. That is the distribution of a Gaussian isotropic field
is invariant under rotation. The definition of isotropy given above is a
generalization of this property for non-Gaussian case. In general, the
isotropy follows and followed by that all higher order moments are also
invariant under rotation. Let us consider homogeneous and isotropic stochastic
field $X\left(  \underline{x}\right)  =X\left(  r,\varphi\right)  $, ($r>0$,
$\varphi\in\left[  0,2\pi\right)  $) on the plane defined by
(\ref{Repr_spherical})
\[
X\left(  r,\varphi\right)  =\sum_{\ell=-\infty}^{\infty}e^{i\ell\varphi}%
\int_{0}^{\infty}J_{\ell}\left(  \rho r\right)  Z_{\ell}\left(  \rho
d\rho\right)  ,
\]
where $Z_{\ell}$ is an array of measures, orthogonal to each other satisfying
(\ref{Z_cov}). In this way an isotropic random field $X\left(  \underline{x}%
\right)  $ can be decomposed into a countable number of mutually uncorrelated
spectral measures defined on the real line instead of on the whole plane,
\cite{adler2010geometr}.

The rotation $g$ takes effect on the 'spherical harmonics' $e^{i\ell m\varphi
}$ ($m=\pm1$), as $\Lambda\left(  g\right)  e^{i\ell m\varphi}=e^{i\ell
m\left(  \varphi+\gamma\right)  }=e^{i\ell m\gamma}e^{i\ell m\varphi}$, since
the $e^{i\ell m\varphi}\ $can be considered as a function of $\varphi$. The
isotropy of $X\left(  r,\varphi\right)  $ implies that the distribution of
$X\left(  r,\varphi\right)  $ does not change under rotations $g\in SO\left(
2\right)  $. Consider
\begin{align*}
\Lambda\left(  g\right)  X\left(  r,\varphi\right)   &  =\sum_{\ell=-\infty
}^{\infty}e^{i\ell\left(  \varphi+\gamma\right)  }\int_{0}^{\infty}J_{\ell
}\left(  \rho r\right)  Z_{\ell}\left(  \rho d\rho\right)  \\
&  =\sum_{\ell=-\infty}^{\infty}e^{i\ell\varphi}\int_{0}^{\infty}J_{\ell
}\left(  \rho r\right)  e^{i\ell\gamma}Z_{\ell}\left(  \rho d\rho\right)  \\
&  =\sum_{\ell=-\infty}^{\infty}e^{i\ell\varphi}\int_{0}^{\infty}J_{\ell
}\left(  \rho r\right)  Z_{\ell}\left(  \rho d\rho\right)  ,
\end{align*}
and because of the above the distribution of $Z_{\ell}\left(  \rho
d\rho\right)  $ and $e^{i\ell\gamma}Z_{\ell}\left(  \rho d\rho\right)  $
should be the same. Now it is evident that for a Gaussian random field
$X\left(  r,\varphi\right)  $, the necessary and sufficient condition of
isotropy is that $Z_{\ell}\left(  \rho d\rho\right)  $ are independent. Indeed
under isotropy assumption we have
\[
\operatorname*{Cum}\left(  Z_{\ell_{1}}\left(  \rho_{1}d\rho_{1}\right)
,Z_{\ell_{2}}\left(  \rho_{2}d\rho_{2}\right)  \right)  =e^{i\left(  \ell
_{1}+\ell_{2}\right)  \gamma}\operatorname*{Cum}\left(  Z_{\ell_{1}}\left(
\rho_{1}d\rho_{1}\right)  ,Z_{\ell_{2}}\left(  \rho_{2}d\rho_{2}\right)
\right)  ,
\]
for each $\gamma$, hence either $\ell_{1}+\ell_{2}=0$, or otherwise
$\operatorname*{Cum}\left(  Z_{\ell_{1}}\left(  \rho_{1}d\rho_{1}\right)
,Z_{\ell_{2}}\left(  \rho_{2}d\rho_{2}\right)  \right)  =0$, and therefore
\[
\operatorname*{Cov}\left(  Z_{\ell_{1}}\left(  \rho_{1}d\rho_{1}\right)
,Z_{\ell_{2}}\left(  \rho_{2}d\rho_{2}\right)  \right)  =0,
\]
unless $\ell_{1}=\ell_{2}$.

In general, we have that under assumption of isotropy the $p^{th}$ order
cumulants satisfy the following equation
\[
\operatorname*{Cum}\left(  Z_{\ell_{1}}\left(  \rho_{1}d\rho_{1}\right)
,\ldots,Z_{\ell_{p}}\left(  \rho_{p}d\rho_{p}\right)  \right)  =e^{i\left(
\ell_{1}+\ell_{2}\cdots+\ell_{p}\right)  \gamma}\operatorname*{Cum}\left(
Z_{\ell_{1}}\left(  \rho_{1}d\rho_{1}\right)  ,\ldots,Z_{\ell_{p}}\left(
\rho_{p}d\rho_{p}\right)  \right)  ,
\]
that is either $\ell_{1}+\ell_{2}\cdots+\ell_{p}=0$, or $\operatorname*{Cum}%
\left(  Z_{\ell_{1}}\left(  \rho_{1}d\rho_{1}\right)  ,Z_{\ell_{2}}\left(
\rho_{2}d\rho_{2}\right)  ,\ldots,Z_{\ell_{p}}\left(  \rho_{p}d\rho
_{p}\right)  \right)  =0$. In turn, if this assumption is satisfied then the
cumulants $\operatorname*{Cum}\left(  Z_{\ell_{1}}\left(  \rho_{1}d\rho
_{1}\right)  ,\ldots,Z_{\ell_{p}}\left(  \rho_{p}d\rho_{p}\right)  \right)  $
are invariant under rotation and the field is isotropic.

\subsection{Spectrum}

In this section we briefly review results already known for second order
spectra of the field $X\left(  r,\varphi\right)  $ before we determine similar
results for bispectra. For notational convenience let us denote the integral
$\int_{0}^{\infty}J_{\ell}\left(  \rho r\right)  Z_{\ell}\left(  \rho
d\rho\right)  $ by $z_{\ell}\left(  r\right)  $, and consider the covariance
\begin{align*}
\operatorname*{Cov}\left(  X\left(  \underline{x}\right)  ,X\left(
\underline{y}\right)  \right)   &  =\operatorname*{Cov}\left(  \sum_{\ell
_{1}=-\infty}^{\ell_{1}=\infty}e^{i\ell_{1}\varphi_{1}}z_{\ell_{1}}\left(
r_{1}\right)  ,\sum_{\ell_{2}=-\infty}^{\ell_{2}=\infty}e^{i\ell_{2}%
\varphi_{2}}z_{\ell_{2}}\left(  r_{2}\right)  \right)  \\
&  =2\pi\int_{0}^{\infty}\sum_{\ell=-\infty}^{\infty}e^{i\ell\left(
\varphi_{1}-\varphi_{2}\right)  }J_{\ell}\left(  \rho r_{1}\right)  J_{\ell
}\left(  \rho r_{2}\right)  F\left(  \rho d\rho\right)  \\
&  =2\pi\int_{0}^{\infty}J_{0}\left(  \rho r\right)  F\left(  \rho
d\rho\right)  ,
\end{align*}
where $r_{1}=\left\vert \underline{x}\right\vert $, $r_{2}=\left\vert
\underline{y}\right\vert $ and $r=\left\vert \underline{x}-\underline{y}%
\right\vert $. In arriving at the above we used the addition formula
\[
J_{0}\left(  \rho r\right)  =\sum_{\ell=-\infty}^{\infty}e^{i\ell\left(
\varphi_{1}-\varphi_{2}\right)  }J_{\ell}\left(  \rho r_{1}\right)  J_{\ell
}\left(  \rho r_{2}\right)  ,
\]
of Bessel functions, see \cite{Erdelyi18c} Tom2, Ch7, 7.6.2.(6),
\cite{Yadrenko1983}. Now one may derive the same result using the properties
of homogeneity and isotropy. 

We are going to apply some special cases of the series expansion given by
(\ref{Repr_spherical}),  namely if the location $\left(  r,\varphi\right)  $
is on the $y$-axis i.e., it points on the direction of the 'North pole'
($N=\left(  0,1\right)  $),
\begin{equation}
X\left(  rN\right)  =\sum_{\ell=-\infty}^{\infty}i^{\ell}\int_{0}^{\infty
}J_{\ell}\left(  \rho r\right)  Z_{\ell}\left(  \rho d\rho\right)
,\label{X_north}%
\end{equation}
and at the origin%
\begin{align}
X\left(  \underline{0}\right)   &  =\int_{\mathbb{R}^{2}}Z\left(
d\underline{\omega}\right)  \nonumber\\
&  =\int_{0}^{\infty}Z_{0}\left(  \rho d\rho\right)  .\label{X_0}%
\end{align}
Let $r=\left\vert \underline{x}-\underline{y}\right\vert $, $\mathcal{C}%
_{2}\left(  r\right)  =\operatorname*{Cov}\left(  X\left(  \underline{x}%
\right)  ,X\left(  \underline{y}\right)  \right)  $, and use the invariance
under translation and rotation to obtain
\begin{align*}
\operatorname*{Cov}\left(  X\left(  \underline{x}\right)  ,X\left(
\underline{y}\right)  \right)   &  =\operatorname*{Cov}\left(  X\left(
\underline{x}-\underline{y}\right)  ,X\left(  0\right)  \right)  \\
&  =\operatorname*{Cov}\left(  X\left(  rN\right)  ,X\left(  0\right)
\right)  \\
&  =2\pi\int_{0}^{\infty}J_{0}\left(  \rho r\right)  F\left(  \rho
d\rho\right)  .
\end{align*}
The above shows one to one correspondence between the second order covariance
and its spectral density function, see \cite{Yadrenko1983}. In particular for
absolutely continuos spectral measure $F\left(  \rho d\rho\right)  =\sigma
^{2}\left\vert A\left(  \rho\right)  \right\vert ^{2}\rho d\rho$ we have
\[
\mathcal{C}_{2}\left(  r\right)  =2\pi\int_{0}^{\infty}J_{0}\left(  \rho
r\right)  \sigma^{2}\left\vert A\left(  \rho\right)  \right\vert ^{2}\rho
d\rho,
\]
in turn
\[
\sigma^{2}\left\vert A\left(  \rho\right)  \right\vert ^{2}=\frac{1}{2\pi}%
\int_{0}^{\infty}J_{0}\left(  \rho r\right)  \mathcal{C}_{2}\left(  r\right)
rdr,
\]
when both integrals exist, see \cite{Brilling74}, \cite{Yaglom-book-87}. The
above property of Hankel transform used above is based on the following
property of Bessel functions
\begin{equation}
\int\limits_{0}^{\infty}J_{\ell}\left(  \rho r\right)  J_{\ell}\left(  \kappa
r\right)  rdr=\frac{\delta\left(  \rho-\kappa\right)  }{\rho}%
,\label{Bessel_Int}%
\end{equation}
where $\delta\left(  \rho-\kappa\right)  $ denotes the Dirac 'function', more
precisely $\delta\left(  \cdot\right)  $ is a distribution (measure), see
\cite{Arfken2001} Sect 11.

\section{Bispectrum}

If the field is not Gaussian then the second order properties do not
characterize the distribution. The next characteristics are, in a row, the
third order moments. The third order structure of a homogeneous and isotropic
stochastic field $X\left(  \underline{x}\right)  $ is described by either the
third order covariances (third order cumulants) in spatial domain or the
bispectrum in frequency domain. Using the spectral representation
(\ref{X_hom}) of $X\left(  \underline{x}\right)  $, we obtain the third order
cumulants (central moments) and it is given by
\begin{align*}
\operatorname*{Cum}\left(  X\left(  \underline{x}_{1}\right)  ,X\left(
\underline{x}_{2}\right)  ,X\left(  \underline{x}_{3}\right)  \right)   &
=\iiint\limits_{\mathbb{R}^{2}\times\mathbb{R}^{2}\times\mathbb{R}^{2}%
}e^{i\left(  \Sigma_{1}^{3}\underline{x}_{k}\cdot\underline{\omega}%
_{k}\right)  }\operatorname*{Cum}\left(  Z\left(  d\underline{\omega}%
_{1}\right)  ,Z\left(  d\underline{\omega}_{2}\right)  ,Z\left(
d\underline{\omega}_{3}\right)  \right) \\
&  =\iiint\limits_{\mathbb{R}^{2}\times\mathbb{R}^{2}\times\mathbb{R}^{2}%
}e^{i\left(  \Sigma_{1}^{3}\underline{x}_{k}\cdot\underline{\omega}%
_{k}\right)  }S_{3}\left(  \underline{\omega}_{1},\underline{\omega}%
_{2},\underline{\omega}_{3}\right)  \delta\left(  \Sigma_{1}^{3}%
\underline{\omega}_{k}\right)
{\textstyle\prod\limits_{k=1}^{3}}
d\underline{\omega}_{k},
\end{align*}
where $S_{3}\left(  \underline{\omega}_{1},\underline{\omega}_{2}%
,\underline{\omega}_{3}\right)  =S_{3}\left(  \underline{\omega}%
_{1},\underline{\omega}_{2},-\underline{\omega}_{1}-\underline{\omega}%
_{2}\right)  $ denotes the bispectral density. Under isotropy for each $g\in
SO\left(  2\right)  $
\begin{align*}
\operatorname*{Cum}\left(  X\left(  g\underline{x}_{1}\right)  ,X\left(
g\underline{x}_{2}\right)  ,X\left(  g\underline{x}_{3}\right)  \right)   &
=\iiint\limits_{\mathbb{R}^{2}\times\mathbb{R}^{2}\times\mathbb{R}^{2}%
}e^{i\left(  \Sigma_{1}^{3}\underline{x}_{k}\cdot\underline{\omega}%
_{k}\right)  }S_{3}\left(  g\underline{\omega}_{1},g\underline{\omega}%
_{2},g\underline{\omega}_{3}\right)  \delta\left(  \Sigma_{1}^{3}%
\underline{\omega}_{k}\right)
{\textstyle\prod\limits_{k=1}^{3}}
d\underline{\omega}_{k}\\
&  =\operatorname*{Cum}\left(  X\left(  \underline{x}_{1}\right)  ,X\left(
\underline{x}_{2}\right)  ,X\left(  \underline{x}_{3}\right)  \right)  ,
\end{align*}
hence $S_{3}\left(  g\underline{\omega}_{1},g\underline{\omega}_{2}%
,g\underline{\omega}_{3}\right)  =S_{3}\left(  \underline{\omega}%
_{1},\underline{\omega}_{2},\underline{\omega}_{3}\right)  =S_{3}\left(
\rho_{1},\rho_{2},\rho_{3}\right)  $. Now we apply the invariance of the third
order covariance
\begin{align*}
\operatorname*{Cum}\left(  X\left(  \underline{x}_{1}\right)  ,X\left(
\underline{x}_{2}\right)  ,X\left(  \underline{x}_{3}\right)  \right)   &
=\operatorname*{Cum}\left(  X\left(  0\right)  ,X\left(  \underline{x}%
_{2}-\underline{x}_{1}\right)  ,X\left(  \underline{x}_{3}-\underline{x}%
_{1}\right)  \right) \\
&  =\operatorname*{Cum}\left(  X\left(  0\right)  ,X\left(  \left\vert
\underline{x}_{2}-\underline{x}_{1}\right\vert N\right)  ,X\left(  g\left(
\underline{x}_{3}-\underline{x}_{1}\right)  \right)  \right)  ,
\end{align*}
where $g$ denotes the rotation carrying $\underline{x}_{2}-\underline{x}_{1}$
into the $\left\vert \underline{x}_{2}-\underline{x}_{1}\right\vert $ times
'North pole' ($N=\left(  0,1\right)  $).

The third order covariance $\operatorname*{Cum}\left(  X\left(  0\right)
,X\left(  \underline{x}_{2}\right)  ,X\left(  \underline{x}_{3}\right)
\right)  $ depends on the length of vectors $\underline{x}_{2}$,
$\underline{x}_{3}$ and the angle $\varphi$ between them, this way a triangle
with vertices $0$, $\underline{x}_{2}$, $\underline{x}_{3}$ is formed with
length of the third side $r_{1}$, such that $r_{1}^{2}=r_{2}^{2}+r_{3}%
^{2}-2r_{2}r_{3}\cos\varphi$. According to this definition of $r_{1}$, we
introduce the notation
\[
\mathcal{C}\left(  r_{1},r_{2},r_{3}\right)  =\operatorname*{Cum}\left(
X\left(  0\right)  ,X\left(  \underline{x}_{2}\right)  ,X\left(
\underline{x}_{3}\right)  \right)  .
\]
We show that the bispectrum $S_{3}\left(  \underline{\omega}_{1}%
,\underline{\omega}_{2},\underline{\omega}_{3}\right)  $ for a homogeneous and
isotropic stochastic field $X\left(  \underline{x}\right)  $ depend on wave
numbers $\rho_{1},\rho_{2},\rho_{3}$, ($\rho_{k}=\left\vert \underline{\omega
}_{k}\right\vert $) only, such that the wave numbers $\rho_{1}$, $\rho_{2}$,
and $\rho_{3}$ satisfy the triangle relation. The angle between sides
$\rho_{2}$ and $\rho_{3}$, will be denoted by $\eta$.

The following theorem shows that the usual connection between the third order
covariances and spectra is valid in a particular form for the third order
covariance and the bispectra as well.

\begin{theorem}
\label{Theo_"D_cov_BISP} The third order covariance $\mathcal{C}\left(
r_{1},r_{2},r_{3}\right)  $ and the corresponding  bispectrum $S_{3}\left(
\rho_{1},\rho_{2},\rho_{3}\right)  $ are given by
\[
\mathcal{C}\left(  r_{1},r_{2},r_{3}\right)  =2\pi\iint\limits_{0}^{\infty
}\int_{0}^{\pi}\left(  J_{0}\left(  w_{+}\right)  +J_{0}\left(  w_{-}\right)
\right)  S_{3}\left(  \rho_{1},\rho_{2},\rho_{3}\right)  d\eta%
{\textstyle\prod\limits_{k=2}^{3}}
\rho_{k}d\rho_{k},
\]
in turn%
\[
S_{3}\left(  \rho_{1},\rho_{2},\rho_{3}\right)  =\frac{1}{\left(  2\pi\right)
^{3}}\iint\limits_{0}^{\infty}\int_{0}^{\pi}\left(  J_{0}\left(  w_{+}\right)
+J_{0}\left(  w_{-}\right)  \right)  \mathcal{C}\left(  r_{1},r_{2}%
,r_{3}\right)  d\varphi%
{\textstyle\prod\limits_{k=2}^{3}}
r_{k}dr_{k},
\]
where $r_{1}^{2}=r_{2}^{2}+r_{3}^{2}-2r_{2}r_{3}\cos\varphi$, $\rho_{1}%
^{2}=\rho_{2}^{2}+\rho_{3}^{2}-2\rho_{2}\rho_{3}\cos\eta$, and%
\begin{align*}
w_{+} &  =\sqrt{\left(  \rho_{2}r_{2}\right)  ^{2}+\left(  \rho_{3}%
r_{3}\right)  ^{2}-2\rho_{2}r_{2}\rho_{3}r_{3}\cos\left(  \left(  \eta
+\varphi\right)  \right)  },\\
w_{-} &  =\sqrt{\left(  \rho_{2}r_{2}\right)  ^{2}+\left(  \rho_{3}%
r_{3}\right)  ^{2}-2\rho_{2}r_{2}\rho_{3}r_{3}\cos\left(  \left(  \eta
-\varphi\right)  \right)  }.
\end{align*}
In the above, we assume that both integrals exist.
\end{theorem}

We note that in some cases it is more convenient to use the transformation
\[
\mathcal{T}\left(  \left.  \eta,\rho_{2},\rho_{3}\right\vert \varphi
,r_{2},r_{3}\right)  =J_{0}\left(  \rho_{2}r_{2}\right)  J_{0}\left(  \rho
_{3}r_{3}\right)  +2\sum_{\ell=1}^{\infty}\cos\left(  \ell\varphi\right)
J_{\ell}\left(  \rho_{2}r_{2}\right)  J_{\ell}\left(  \rho_{3}r_{3}\right)
\cos\left(  \ell\eta\right)  ,
\]
between the bispectrum $S_{3}\left(  \rho_{1},\rho_{2},\rho_{3}\right)  $ and
the third order covariance $\mathcal{C}_{3}\left(  r_{1},r_{2},r_{3}\right)
$, namely.
\begin{align}
\mathcal{C}\left(  r_{1},r_{2},r_{3}\right)   &  =4\pi\iint\limits_{0}%
^{\infty}\int_{0}^{\pi}\mathcal{T}\left(  \left.  \eta,\rho_{2},\rho
_{3}\right\vert \varphi,r_{2},r_{3}\right)  S_{3}\left(  \rho_{1},\rho
_{2},\rho_{3}\right)  d\eta%
{\textstyle\prod\limits_{k=2}^{3}}
\rho_{k}d\rho_{k},\label{Bicov_T_transf}\\
S_{3}\left(  \rho_{1},\rho_{2},\rho_{3}\right)   &  =\frac{1}{4\pi^{3}}%
\iint\limits_{0}^{\infty}\int_{0}^{\pi}\mathcal{T}\left(  \left.  \eta
,\rho_{2},\rho_{3}\right\vert \varphi,r_{2},r_{3}\right)  \mathcal{C}\left(
r_{1},r_{2},r_{3}\right)  d\varphi%
{\textstyle\prod\limits_{k=2}^{3}}
r_{k}dr_{k}.\label{Bisp_T_transf}%
\end{align}

\begin{proof}
We use the particular representations (\ref{X_0}), (\ref{X_north}) and obtain
\begin{multline*}
\operatorname*{Cum}\left(  X\left(  0\right)  ,X\left(  r_{2}N\right)
,X\left(  \underline{x}_{3}\right)  \right)  \\
=\sum_{\ell_{2},\ell_{3}=-\infty}^{\infty}\iiint\limits_{0}^{\infty}%
i^{\ell_{2}}e^{i\ell_{3}\varphi_{3}}J_{\ell_{2}}\left(  \rho_{2}r_{2}\right)
J_{\ell_{3}}\left(  \rho_{3}r_{3}\right)  \operatorname*{Cum}\left(
Z_{0}\left(  \rho_{1}d\rho_{1}\right)  ,Z_{\ell_{2}}\left(  \rho_{2}d\rho
_{2}\right)  ,Z_{\ell_{3}}\left(  \rho_{3}d\rho_{3}\right)  \right)  \\
=\sum_{\ell=-\infty}^{\infty}\iiint\limits_{0}^{\infty}e^{i\ell\varphi}%
J_{\ell}\left(  \rho_{2}r_{2}\right)  J_{-\ell}\left(  \rho_{3}r_{3}\right)
\operatorname*{Cum}\left(  Z_{0}\left(  \rho_{1}d\rho_{1}\right)  ,Z_{\ell
}\left(  \rho_{2}d\rho_{2}\right)  ,Z_{-\ell}\left(  \rho_{3}d\rho_{3}\right)
\right)
\end{multline*}
where $\varphi=\pi/2-\varphi_{3}$, is the angle between $N$ and $\underline{x}%
_{3}$. The third order cumulant of the stochastic spectral measure $Z\left(
d\underline{\omega}\right)  $ of the homogeneous field $X\left(
\underline{x}\right)  $ is given by
\begin{align*}
\operatorname*{Cum}\left(  Z\left(  d\underline{\omega}_{1}\right)  ,Z\left(
d\underline{\omega}_{2}\right)  ,Z\left(  d\underline{\omega}_{3}\right)
\right)   &  =\delta\left(  \Sigma_{1}^{3}\underline{\omega}_{k}\right)
S_{3}\left(  \underline{\omega}_{1},\underline{\omega}_{2},\underline{\omega
}_{3}\right)  d\underline{\omega}_{1}d\underline{\omega}_{2}d\underline{\omega
}_{3}\\
&  =\delta\left(  \Sigma_{1}^{3}\rho_{k}\underline{\widehat{\omega}}%
_{k}\right)  S_{3}\left(  \rho_{1},\rho_{2},\rho_{3}\right)
{\textstyle\prod\limits_{k=1}^{3}}
\Omega\left(  d\underline{\widehat{\omega}}_{k}\right)  \rho_{k}d\rho_{k},
\end{align*}
where $\underline{\widehat{\omega}}_{k}=\underline{\omega}_{k}/\left\vert
\underline{\omega}_{k}\right\vert $. The stochastic spectral measures
$Z_{\ell}\left(  \rho d\rho\right)  $ are related to $Z\left(
d\underline{\omega}\right)  $ by
\[
Z_{\ell}\left(  \rho d\rho\right)  =\int_{0}^{2\pi}i^{\ell}e^{-i\ell\eta
\ }Z\left(  d\eta\rho d\rho\right)  ,
\]
therefore
\begin{align}
\operatorname*{Cum}\left(  Z_{0}\left(  \rho_{1}d\rho_{1}\right)  ,Z_{\ell
}\left(  \rho_{2}d\rho_{2}\right)  ,Z_{-\ell}\left(  \rho_{3}d\rho_{3}\right)
\right)   &  =S_{3}\left(  \rho_{1},\rho_{2},\rho_{3}\right)  \label{Cum3_Z_L}%
\\
&  \times\iiint\limits_{0}^{2\pi}e^{-i\ell\left(  \eta_{3}-\eta_{2}\right)
\ }\delta\left(  \Sigma_{1}^{3}\rho_{k}\underline{\widehat{\omega}}%
_{k}\right)
{\textstyle\prod\limits_{k=1}^{3}}
\rho_{k}d\rho_{k}d\eta_{k}.\nonumber
\end{align}
In order to understand the usefulness of the Dirac 'function' in polar
coordinates we express it by an integral through the Jacobi-Anger expansion
(\ref{Expans_Jacobi_Anger}). Since the Dirac 'function' is a measure we apply
here the theory of generalized functions to  obtain
\begin{align}
\delta\left(  \Sigma_{1}^{3}\rho_{k}\underline{\widehat{\omega}}_{k}\right)
&  =\frac{1}{\left(  2\pi\right)  ^{2}}\iint\limits_{\mathbb{R}^{2}%
}e^{i\left(  \underline{\lambda}\cdot\Sigma_{1}^{3}\underline{\omega}%
_{k}\right)  }d\underline{\lambda}\label{Dirac3}\\
&  =\frac{1}{\left(  2\pi\right)  ^{2}}\int_{0}^{\infty}\int_{0}^{2\pi}%
{\displaystyle\prod\limits_{k=1}^{3}}
\sum_{\ell_{k}=-\infty}^{\infty}i^{\ell_{k}}J_{\ell_{k}}\left(  \rho
_{k}\lambda\right)  e^{i\ell_{k}\left(  \eta_{k}-\xi\right)  }\lambda d\lambda
d\xi\nonumber\\
&  =\frac{1}{2\pi}\int_{0}^{\infty}\sum_{\ell_{1},\ell_{2}=-\infty}^{\infty
}e^{i\left(  \ell_{1}\left(  \eta_{1}-\eta_{3}\right)  +\ell_{2}\left(
\eta_{2}-\eta_{3}\right)  \right)  }J_{\ell_{1}}\left(  \rho_{1}%
\lambda\right)  J_{\ell_{2}}\left(  \rho_{2}\lambda\right)  J_{-\ell_{1}%
-\ell_{2}}\left(  \rho_{3}\lambda\right)  \lambda d\lambda.\nonumber
\end{align}
Now substitute (\ref{Dirac3}) into (\ref{Cum3_Z_L}). Because
\[
\iiint\limits_{0}^{2\pi}e^{-i\ell\left(  \eta_{3}-\eta_{2}\right)  \ }%
{\textstyle\prod\limits_{k=1}^{3}}
e^{i\ell_{k}\left(  \eta_{k}-\xi\right)  }d\eta_{k}=\delta_{\ell_{1}}%
\delta_{\ell_{2}+\ell}\delta_{\ell_{3}-\ell}\left(  2\pi\right)  ^{3},
\]
we get
\[
\iiint\limits_{0}^{2\pi}e^{-i\ell\left(  \eta_{3}-\eta_{2}\right)  \ }%
\delta\left(  \Sigma_{1}^{3}\rho_{k}\underline{\widehat{\omega}}_{k}\right)
d\eta_{1}d\eta_{2}d\eta_{3}=\left(  2\pi\right)  ^{2}\int_{0}^{\infty}%
J_{0}\left(  \rho_{1}\lambda\right)  J_{-\ell}\left(  \rho_{2}\lambda\right)
J_{\ell}\left(  \rho_{3}\lambda\right)  \lambda d\lambda.
\]
This integral can be evaluated, if $\left\vert \rho_{2}-\rho_{3}\right\vert
<\rho_{1}<\rho_{2}+\rho_{3}$, and let us denote $R=\left(  \rho_{2}^{2}%
+\rho_{3}^{2}-\rho_{1}^{2}\right)  /\left(  2\rho_{2}\rho_{3}\right)  $,
$\rho_{1}^{2}=\rho_{2}^{2}+\rho_{3}^{2}-2\rho_{2}\rho_{3}\cos\left(
\eta\right)  $, then
\[
\int_{0}^{\infty}J_{0}\left(  \rho_{1}\lambda\right)  J_{\ell}\left(  \rho
_{2}\lambda\right)  J_{\ell}\left(  \rho_{3}\lambda\right)  \lambda
d\lambda=\frac{\cos\left(  \ell\arccos\left(  R\right)  \right)  }{\pi\rho
_{2}\rho_{3}\sqrt{1-R^{2}}},
\]
see \cite{Prudnikov14} Tom. II, 2.12.41.16, hence
\[
\iiint\limits_{0}^{2\pi}e^{-i\ell\left(  \eta_{1}-\eta_{2}\right)  \ }%
\delta\left(  \Sigma_{1}^{3}\rho_{k}\underline{\widehat{\omega}}_{k}\right)
d\eta_{1}d\eta_{2}d\eta_{3}=4\pi\left(  -1\right)  ^{\ell}\frac{\cos\left(
\ell\arccos\left(  R\right)  \right)  }{\rho_{2}\rho_{3}\sqrt{1-R^{2}}}.
\]
Since $J_{\ell}=\left(  -1\right)  ^{\ell}J_{-\ell}$ we have
\begin{multline*}
\operatorname*{Cum}\left(  Z_{0}\left(  \rho_{1}d\rho_{1}\right)  ,Z_{\ell
}\left(  \rho_{2}d\rho_{2}\right)  ,Z_{-\ell}\left(  \rho_{3}d\rho_{3}\right)
\right)  \\
=4\pi\left(  -1\right)  ^{\ell}\delta\left(  \rho\triangle\right)  \frac
{\cos\left(  \ell\arccos\left(  R\right)  \right)  }{\rho_{2}\rho_{3}%
\sqrt{1-R^{2}}}S_{3}\left(  \rho_{1},\rho_{2},\rho_{3}\right)
{\textstyle\prod\limits_{k=1}^{3}}
\rho_{k}d\rho_{k},
\end{multline*}
where $\delta\left(  \rho\triangle\right)  =\delta\left(  \rho_{2}^{2}%
+\rho_{3}^{2}-2\rho_{2}\rho_{3}\cos\eta-\rho_{1}^{2}\right)  $ and therefore
the wave numbers $\rho_{1}$, $\rho_{2}$, and $\rho_{3}$ should satisfy the
triangle relation. Now can obtain the following expression for the third order
covariance
\begin{align}
&  \operatorname*{Cum}\left(  X\left(  0\right)  ,X\left(  r_{2}N\right)
,X\left(  \underline{x}_{3}\right)  \right)  \nonumber\\
&  =4\pi\sum_{\ell=-\infty}^{\infty}\iiint\limits_{0}^{\infty}e^{-i\ell
\varphi}J_{\ell}\left(  \rho_{2}r_{2}\right)  J_{-\ell}\left(  \rho_{3}%
r_{3}\right)  \left(  -1\right)  ^{\ell}\delta\left(  \rho\triangle\right)
\frac{\cos\left(  \ell\arccos\left(  R\right)  \right)  }{\rho_{2}\rho
_{3}\sqrt{1-R^{2}}}S_{3}\left(  \rho_{1},\rho_{2},\rho_{3}\right)
{\textstyle\prod\limits_{k=1}^{3}}
\rho_{k}d\rho_{k}\nonumber\\
&  =4\pi\sum_{\ell=-\infty}^{\infty}\iint\limits_{0}^{\infty}\int_{\left\vert
\rho_{2}-\rho_{3}\right\vert }^{\rho_{2}+\rho_{3}}e^{-i\ell\varphi}J_{\ell
}\left(  \rho_{2}r_{2}\right)  J_{\ell}\left(  \rho_{3}r_{3}\right)
\frac{\cos\left(  \ell\arccos\left(  R\right)  \right)  }{\rho_{2}\rho
_{3}\sqrt{1-R^{2}}}S_{3}\left(  \rho_{1},\rho_{2},\rho_{3}\right)
{\textstyle\prod\limits_{k=1}^{3}}
\rho_{k}d\rho_{k}\nonumber\\
&  =4\pi\sum_{\ell=-\infty}^{\infty}\iint\limits_{0}^{\infty}\int_{0}^{\pi
}e^{-i\ell\varphi}J_{\ell}\left(  \rho_{2}r_{2}\right)  J_{\ell}\left(
\rho_{3}r_{3}\right)  \frac{\cos\left(  \ell\eta\right)  }{\sqrt{1-\cos
^{2}\left(  \eta\right)  }}S_{3}\left(  \rho_{1},\rho_{2},\rho_{3}\right)
\sin\eta d\eta%
{\textstyle\prod\limits_{k=2}^{3}}
\rho_{k}d\rho_{k}\nonumber\\
&  =4\pi\iint\limits_{0}^{\infty}\int_{0}^{\pi}\sum_{\ell=-\infty}^{\infty
}e^{-i\ell\varphi}J_{\ell}\left(  \rho_{2}r_{2}\right)  J_{\ell}\left(
\rho_{3}r_{3}\right)  \cos\left(  \ell\eta\right)  S_{3}\left(  \rho_{1}%
,\rho_{2},\rho_{3}\right)  d\eta%
{\textstyle\prod\limits_{k=2}^{3}}
\rho_{k}d\rho_{k},\label{Transf_T_s}%
\end{align}
where $\rho_{1}^{2}=\rho_{2}^{2}+\rho_{3}^{2}-2\rho_{2}\rho_{3}\cos\left(
\eta\right)  $ and
\[
\frac{\rho_{1}d\rho_{1}}{d\eta}=\rho_{2}\rho_{3}\sin\eta.
\]
The function in the last row of the expression (\ref{Transf_T_s}) is
\begin{align}
\mathcal{T}\left(  \left.  \eta,\rho_{2},\rho_{3}\right\vert \varphi
,r_{2},r_{3}\right)   &  =\sum_{\ell=-\infty}^{\infty}e^{-i\ell\varphi}%
J_{\ell}\left(  \rho_{2}r_{2}\right)  J_{\ell}\left(  \rho_{3}r_{3}\right)
\cos\left(  \ell\eta\right)  \nonumber\\
&  =\sum_{\ell=-\infty}^{\infty}\cos\left(  \ell\varphi\right)  J_{\ell
}\left(  \rho_{2}r_{2}\right)  J_{\ell}\left(  \rho_{3}r_{3}\right)
\cos\left(  \ell\eta\right)  \nonumber\\
&  =J_{0}\left(  \rho_{2}r_{2}\right)  J_{0}\left(  \rho_{3}r_{3}\right)
+2\sum_{\ell=1}^{\infty}\cos\left(  \ell\varphi\right)  J_{\ell}\left(
\rho_{2}r_{2}\right)  J_{\ell}\left(  \rho_{3}r_{3}\right)  \cos\left(
\ell\eta\right)  .\label{Transf_T}%
\end{align}
The above gives the transformation of the bispectrum $S_{3}\left(  \rho
_{1},\rho_{2},\rho_{3}\right)  $ from the third order covariance
$\mathcal{C}_{3}\left(  r_{1},r_{2},r_{3}\right)  $. Note that both angles
$\varphi$ and $\eta$ together with two sides define the third side $\rho_{1}$
and $r_{1}$ of the triangles, given by wave numbers $\left(  \rho_{1},\rho
_{2},\rho_{3}\right)  $ and distances $\left(  r_{1},r_{2},r_{3}\right)  $.
The transformation $\mathcal{T}$ can be simplified by using
\begin{multline*}
\mathcal{T}\left(  \left.  \eta,\rho_{2},\rho_{3}\right\vert \varphi
,r_{2},r_{3}\right)  =\frac{1}{2}\left(  J_{0}\left(  \rho_{2}r_{2}\right)
J_{0}\left(  \rho_{3}r_{3}\right)  +2\sum_{\ell=1}^{\infty}J_{\ell}\left(
\rho_{2}r_{2}\right)  J_{\ell}\left(  \rho_{3}r_{3}\right)  \cos\left(
\ell\left(  \eta+\varphi\right)  \right)  \right)  \\
+\frac{1}{2}\left(  J_{0}\left(  \rho_{2}r_{2}\right)  J_{0}\left(  \rho
_{3}r_{3}\right)  +2\sum_{\ell=1}^{\infty}J_{\ell}\left(  \rho_{2}%
r_{2}\right)  J_{\ell}\left(  \rho_{3}r_{3}\right)  \cos\left(  \ell\left(
\eta-\varphi\right)  \right)  \right)  \\
=\frac{1}{2}\left(  J_{0}\left(  w_{+}\right)  +J_{0}\left(  w_{-}\right)
\right)  ,
\end{multline*}
where $w_{+}=\sqrt{\left(  \rho_{2}r_{2}\right)  ^{2}+\left(  \rho_{3}%
r_{3}\right)  ^{2}-2\rho_{2}r_{2}\rho_{3}r_{3}\cos\left(  \left(  \eta
+\varphi\right)  \right)  }$, \newline$w_{-}=\sqrt{\left(  \rho_{2}%
r_{2}\right)  ^{2}+\left(  \rho_{3}r_{3}\right)  ^{2}-2\rho_{2}r_{2}\rho
_{3}r_{3}\cos\left(  \left(  \eta-\varphi\right)  \right)  }$, see
\cite{Erdelyi18c} Tom.2, Ch7.7.15, for the formula
\begin{align*}
J_{0}\left(  z\right)  J_{0}\left(  u\right)  +2\sum_{\ell=1}^{\infty}J_{\ell
}\left(  z\right)  J_{\ell}\left(  u\right)  \cos\left(  \ell\vartheta\right)
&  =J_{0}\left(  w\right)  ,\\
w &  =\sqrt{z^{2}+u^{2}-2zu\cos\left(  \vartheta\right)  }.
\end{align*}
We have now established a relationship between the third order covariance and
the bispectrum, namely
\[
\operatorname*{Cum}\left(  X\left(  0\right)  ,X\left(  r_{2}N\right)
,X\left(  \underline{x}_{3}\right)  \right)  =2\pi\iint\limits_{0}^{\infty
}\int_{0}^{\pi}\left(  J_{0}\left(  w_{+}\right)  +J_{0}\left(  w_{-}\right)
\right)  S_{3}\left(  \rho_{1},\rho_{2},\rho_{3}\right)  d\eta%
{\textstyle\prod\limits_{k=2}^{3}}
\rho_{k}d\rho_{k}.
\]
Next we show that the inversion formula
\[
S_{3}\left(  \rho_{1},\rho_{2},\rho_{3}\right)  =\frac{1}{\left(  2\pi\right)
^{3}}\iint\limits_{0}^{\infty}\int_{0}^{\pi}\left(  J_{0}\left(  w_{+}\right)
+J_{0}\left(  w_{-}\right)  \right)  \mathcal{C}\left(  r_{1},r_{2}%
,r_{3}\right)  d\varphi%
{\textstyle\prod\limits_{k=2}^{3}}
r_{k}dr_{k},
\]
is also valid. Consider the integral
\begin{align}
I\left(  \left.  \rho_{1},\rho_{2},\rho_{3}\right\vert \kappa_{1},\kappa
_{2},\kappa_{3}\right)   &  =\iint\limits_{0}^{\infty}\int_{0}^{\pi}\left(
J_{0}\left(  \rho_{2}r_{2}\right)  J_{0}\left(  \rho_{3}r_{3}\right)
+2\sum_{\ell=1}^{\infty}\cos\left(  \ell\varphi\right)  J_{\ell}\left(
\rho_{2}r_{2}\right)  J_{\ell}\left(  \rho_{3}r_{3}\right)  \cos\left(
\ell\eta\right)  \right)  \nonumber\\
&  \hspace{-0.9in}\times\left(  J_{0}\left(  \kappa_{2}r_{2}\right)
J_{0}\left(  \kappa_{3}r_{3}\right)  +2\sum_{\ell=1}^{\infty}\cos\left(
\ell\varphi\right)  J_{\ell}\left(  \kappa_{2}r_{2}\right)  J_{\ell}\left(
\kappa_{3}r_{3}\right)  \cos\left(  \ell\vartheta\right)  \right)  d\varphi%
{\textstyle\prod\limits_{k=2}^{3}}
r_{k}dr_{k}.\label{Int_I}%
\end{align}
WE notice  first that $\cos\left(  \ell\varphi\right)  $ is an orthogonal
system on $\left[  0,\pi\right]  $, i.e.
\[
\int_{0}^{\pi}\cos\left(  \ell_{1}\varphi\right)  \cos\left(  \ell_{2}%
\varphi\right)  d\varphi=\delta_{\ell_{1}=\ell_{2}}\left\{
\begin{array}
[c]{ccc}%
\pi & \text{if } & \ell_{1}=0\\
\frac{\pi}{2} & \text{if } & \ell_{1}\neq0
\end{array}
\right.  ,
\]
then we integrate (\ref{Int_I}) with respect to $\varphi$ and obtain
\begin{multline*}
I\left(  \left.  \rho_{1},\rho_{2},\rho_{3}\right\vert \kappa_{1},\kappa
_{2},\kappa_{3}\right)  =\pi\int\limits_{0}^{\infty}J_{0}\left(  \rho_{2}%
r_{2}\right)  J_{0}\left(  \kappa_{2}r_{2}\right)  r_{2}dr_{2}\int%
\limits_{0}^{\infty}J_{0}\left(  \rho_{3}r_{3}\right)  J_{0}\left(  \kappa
_{3}r_{3}\right)  r_{3}dr_{3}\\
+2\pi\sum_{\ell=1}^{\infty}\cos\left(  \ell\eta\right)  \cos\left(
\ell\vartheta\right)  \int\limits_{0}^{\infty}J_{\ell}\left(  \rho_{2}%
r_{2}\right)  J_{\ell}\left(  \kappa_{2}r_{2}\right)  r_{2}dr_{2}%
\int\limits_{0}^{\infty}J_{\ell}\left(  \rho_{3}r_{3}\right)  J_{\ell}\left(
\kappa_{3}r_{3}\right)  r_{3}dr_{3}.
\end{multline*}
Using the integral of Bessel functions (\ref{Bessel_Int}), we can show
\begin{align*}
I\left(  \left.  \rho_{1},\rho_{2},\rho_{3}\right\vert \kappa_{1},\kappa
_{2},\kappa_{3}\right)   &  =\pi^{2}\delta\left(  \rho_{2}-\kappa_{2}\right)
\delta\left(  \rho_{3}-\kappa_{3}\right)  \left(  \frac{1}{\pi}+\frac{2}{\pi
}\sum_{\ell=1}^{\infty}\cos\left(  \ell\eta\right)  \cos\left(  \ell
\vartheta\right)  \right)  \\
&  =\pi^{2}\delta\left(  \rho_{2}-\kappa_{2}\right)  \delta\left(  \rho
_{3}-\kappa_{3}\right)  \delta\left(  \eta-\vartheta\right)  ,
\end{align*}
since $1/\pi$ and $\sqrt{\pi/2}\cos\left(  \ell\eta\right)  $ forms an
orthonormal system, see \cite{Arfken2001} Sect.11.
\end{proof}

We now consider  models to define above random precesses.

\begin{example}
The spatial white noise $\partial W\left(  r,\varphi\right)  $, on the plane
is given as a generalized field by the series representation
\[
\partial W\left(  r,\varphi\right)  =\sum_{\ell=-\infty}^{\infty}%
e^{i\ell\varphi}\int_{0}^{\infty}J_{\ell}\left(  \rho r\right)  W_{\ell
}\left(  \rho d\rho\right)  ,
\]
where $W_{\ell}\left(  \rho d\rho\right)  $ with $E\left\vert W_{\ell}\left(
\rho d\rho\right)  \right\vert ^{2}=\sigma^{2}\rho d\rho$, see
\cite{Yaglom-book-87}, \cite{Yadrenko1983}. We define the Laplacian field on
the plane by the equation
\begin{equation}
\left(  \nabla^{2}-c^{2}\right)  X\left(  r,\varphi\right)  =\partial W\left(
r,\varphi\right)  ,\label{Equ_Laplace}%
\end{equation}
where
\[
\nabla^{2}=\frac{1}{r}\frac{\partial}{\partial r}+\frac{\partial^{2}}{\partial
r^{2}}+\frac{1}{r^{2}}\frac{\partial^{2}}{\partial\varphi^{2}},
\]
is the Laplacian operator in terms of spherical coordinates. Now we have
\begin{align*}
\nabla^{2}e^{i\ell\varphi}\int_{0}^{\infty}J_{\ell}\left(  \rho r\right)
Z_{\ell}\left(  \rho d\rho\right)   &  =e^{i\ell\varphi}\int_{0}^{\infty
}\left[  \frac{d^{2}}{dr^{2}}+\frac{1}{r}\frac{d}{dr}-\frac{\ell^{2}}{r^{2}%
}\right]  J_{\ell}\left(  \rho r\right)  Z_{\ell}\left(  \rho d\rho\right)  \\
&  =-e^{i\ell\varphi}\int_{0}^{\infty}\rho^{2}J_{\ell}\left(  \rho r\right)
Z_{\ell}\left(  \rho d\rho\right)  ,
\end{align*}
hence
\[
\left(  \nabla^{2}-c^{2}\right)  e^{i\ell\varphi}\int_{0}^{\infty}J_{\ell
}\left(  \rho r\right)  Z_{\ell}\left(  \rho d\rho\right)  =-e^{i\ell\varphi
}\int_{0}^{\infty}\left(  \rho^{2}+c^{2}\right)  J_{\ell}\left(  \rho
r\right)  Z_{\ell}\left(  \rho d\rho\right)  .
\]
Let us compare the terms in the integrands  of the integrals of the equation
(\ref{Equ_Laplace}), and obtain
\begin{align*}
-\left(  \rho^{2}+c^{2}\right)  Z_{\ell}\left(  \rho d\rho\right)   &
=W_{\ell}\left(  \rho d\rho\right)  ,\\
\left(  \rho^{2}+c^{2}\right)  ^{2}2\pi F\left(  \rho d\rho\right)   &
=\sigma^{2}\rho d\rho,\\
F\left(  \rho d\rho\right)   &  =\frac{\sigma^{2}\rho d\rho}{2\pi\left(
\rho^{2}+c^{2}\right)  ^{2}}.
\end{align*}
Hence the covariance is obtained by inversion,
\begin{align*}
\operatorname*{Cov}\left(  X\left(  \underline{x}\right)  ,X\left(
\underline{y}\right)  \right)   &  =2\pi\int_{0}^{\infty}J_{0}\left(  \rho
r\right)  F\left(  \rho d\rho\right)  \\
&  =\int_{0}^{\infty}J_{0}\left(  \rho r\right)  \frac{\sigma^{2}\rho d\rho
}{\left(  \rho^{2}+c^{2}\right)  ^{2}}\\
&  =\sigma^{2}\frac{rK_{-1}\left(  cr\right)  }{2c}.
\end{align*}
This covariance belongs to Mat\'{e}rn Class, see \cite{Whittle1954}
($K_{-1}=K_{1})$, in terms of modified Bessel (Hankel) function. The
bispectrum of the process defined by the Laplacian model is given by
\[
S_{3}\left(  \rho_{1},\rho_{2},\rho_{3}\right)  =%
{\textstyle\prod\limits_{k=1}^{3}}
\frac{\sigma^{2}}{\rho_{k}^{2}+c^{2}},
\]
where $\rho_{1}^{2}=\rho_{2}^{2}+\rho_{3}^{2}-2\rho_{2}\rho_{3}\cos\left(
\eta\right)  $. We express the third order covariances according to the
Theorem \ref{Theo_"D_cov_BISP}, applying the transformation (\ref{Transf_T})
to the bispectrum. If $\left(  a/b\right)  ^{2}<1,$ we have
\[
\int_{0}^{\pi}\frac{\cos\left(  \ell\eta\right)  }{b+a\cos\left(  \eta\right)
}d\eta=\frac{\pi}{\sqrt{b^{2}-a^{2}}}\left(  \frac{\sqrt{b^{2}-a^{2}}-b}%
{a}\right)  ^{\ell},
\]
see \cite{GradshteynRyzhik} 3.613.1. Put
\begin{align*}
a &  =-2\rho_{2}\rho_{3},\quad b=\rho_{2}^{2}+\rho_{3}^{2}+c^{2}\\
b^{2}-a^{2} &  =\left(  \left(  \rho_{2}-\rho_{3}\right)  ^{2}+c^{2}\right)
\left(  \left(  \rho_{2}+\rho_{3}\right)  ^{2}+c^{2}\right)  ,
\end{align*}
hence
\[
\int_{0}^{\pi}\frac{\cos\left(  \ell\eta\right)  }{\rho_{2}^{2}+\rho_{3}%
^{2}-2\rho_{2}\rho_{3}\cos\left(  \eta\right)  +c^{2}}d\eta=\frac{\pi}%
{\sqrt{b^{2}-a^{2}}}\left(  \frac{\sqrt{b^{2}-a^{2}}-b}{a}\right)  ^{\ell}.
\]
In this way we arrive at the Fourier expansion of the third order covariances
\[
\mathcal{C}\left(  r_{1},r_{2},r_{3}\right)  =f_{0}+2\sum_{\ell=1}^{\infty
}f_{l}\cos\left(  \ell\varphi\right)  ,
\]
with coefficients
\[
f_{l}=4\pi\iint\limits_{0}^{\infty}\frac{\pi}{\sqrt{b^{2}-a^{2}}}\left(
\frac{\sqrt{b^{2}-a^{2}}-b}{a}\right)  ^{\ell}J_{\ell}\left(  \rho_{2}%
r_{2}\right)  J_{\ell}\left(  \rho_{3}r_{3}\right)
{\textstyle\prod\limits_{k=2}^{3}}
\rho_{k}d\rho_{k}.
\]

\end{example}

\begin{acknowledgement}
The publication was supported by the T\'{A}MOP-4.2.2.C-11/1/KONV-2012-0001
project. The project has been supported by the European Union, co-financed by
the European Social Fund.
\end{acknowledgement}

\bibliographystyle{aalpha}
\bibliography{00BiblMM13}

\end{document}